\documentclass{article}

\usepackage{amsmath,amssymb,amsthm}

\newtheorem{theorem}{Theorem}
\newtheorem{definition}{Definition}
\theoremstyle{remark}
\newtheorem{remark}{Remark}

\providecommand{\ZBL}{\relax\ifhmode\unskip\space\fi \textsf{Zbl} }
\providecommand{\JFM}{\relax\ifhmode\unskip\space\fi \textsf{JFM} }
\providecommand{\MR}{\relax\ifhmode\unskip\space\fi \textsf{MR} }

\begin{document}

\title{A Proper Extension of Noether's Symmetry Theorem
for Nonsmooth Extremals\\ of the Calculus of Variations\footnote{Research
Report CM03/I-11}}

\author{Delfim F. M. Torres\\
        \texttt{delfim@mat.ua.pt}}

\date{Department of Mathematics\\
      University of Aveiro\\
      3810-193 Aveiro, Portugal\\
      \texttt{http://www.mat.ua.pt/delfim}}

\maketitle


\begin{abstract}
For nonsmooth Euler-Lagrange extremals, Noether's
conservation laws cease to be valid.
We show that Emmy Noether's theorem of the calculus
of variations is still valid in the wider class of Lipschitz
functions, as long as one restrict the Euler-Lagrange extremals to those
which satisfy the DuBois-Reymond necessary condition. In
the smooth case all Euler-Lagrange extremals are DuBois-Reymond extremals,
and the result gives a proper extension of the classical Noether's theorem.
This is in contrast with the recent developments of Noether's
symmetry theorems to the optimal control setting, which give rise to
non-proper extensions when specified for the problems
of the calculus of variations.
\end{abstract}


\vspace*{0.5cm}



\vspace*{0.3cm}



\section{Introduction}

Let $L(t,x,v)$ be a given
$C^1\left([a,b]\times\mathbb{R}^n\times\mathbb{R}^n;\mathbb{R}\right)$
function (the Lagrangian). The fundamental problem of the calculus of
variations consists to minimize the integral functional
\begin{equation}
\label{eqT:J}
J\left[x(\cdot)\right] =
\int_a^b L\left(t,x(t),\dot{x}(t)\right) \mathrm{d}t
\end{equation}
over a certain class $\mathcal{X}$ of functions
$x : [a,b] \rightarrow \mathbb{R}^n$ satisfying the boundary
conditions $x(a) = \alpha$, $x(b) = \beta$. The problem is
usually solved with the help of the famous Euler-Lagrange
equation,
\begin{equation}
\label{eqT:EL}
\frac{\mathrm{d}}{\mathrm{d}t}
\frac{\partial L}{\partial v}\left(t,x(t),\dot{x}(t)\right)
= \frac{\partial L}{\partial x}\left(t,x(t),\dot{x}(t)\right) \, ,
\end{equation}
which is a first-order necessary optimality condition.
Each solution of \eqref{eqT:EL} is called an Euler-Lagrange extremal.
Condition \eqref{eqT:EL} is obtained in most textbooks from the
assumption that minimizers are smooth, or assuming they are piecewise
smooth functions. In this last situation the Euler-Lagrange equations
are interpreted as holding everywhere except possibly at finitely many points.

In 1918 Emmy Noether \cite{JFM46.0770.01,MR53:10538}
established a general theorem asserting that the
invariance of the integral functional \eqref{eqT:J} under a
group of transformations depending smoothly on a parameter $s$,
implies the existence of a conserved quantity along the Euler-Lagrange
extremals. As corollaries, all the conservation laws known to classical mechanics
are easily obtained. For a survey of Noether's theorem and its
generalizations see \cite{MR83c:70020}.
Noether's theorem, as is found in the many literature
of physics, calculus of variations and optimal control,
is formulated with $\mathcal{X}$ being smooth. A typical example
is $x(\cdot) \in \mathcal{X} = C^2$
(\textrm{cf. e.g.} \cite{MR89e:49002,MR28:3353,MR2000m:49002}).

Given that the Euler-Lagrange
equation \eqref{eqT:EL} makes sense when $x(\cdot)$ has merely
essentially bounded derivative -- the biggest class $\mathcal{X}$ for which
\eqref{eqT:EL} is still valid is the class $Lip$ of Lipschitz
functions (\textrm{cf. e.g.} \cite[\S 2.2]{MR85c:49001})\footnote{In this situation
equations \eqref{eqT:EL} are interpreted in the \emph{almost everywhere} sense.} -- it
is expected that the conclusion of Noether's theorem can still be
defended in the wider class of Lipschitz functions. This is
indeed the case, as it follows from the Pontryagin
maximum principle and the results in \cite{MR1901565,delfimEJC}.
As far as for the fundamental problem of the calculus of variations
the Pontryagin maximum principle reduces to the Euler-Lagrange
necessary condition \eqref{eqT:EL} and to the
Weierstrass necessary condition
\begin{equation}
\label{eqT:Weiers}
L\left(t,x(t),v\right) - L\left(t,x(t),\dot{x}(t)\right) \ge
\frac{\partial L}{\partial v}\left(t,x(t),\dot{x}(t)\right)
\cdot \left(v - \dot{x}(t)\right) \quad \forall \
v \in \mathbb{R}^n \, ,
\end{equation}
which are distinct necessary
conditions even in the $C^2$-smooth case, this does not give
a proper extension of Noether's theorem to the class of Lipschitz functions
(the generalization does not reduce to the classical formulation when
$\mathcal{X} = C^2$, since we are restricting the set of Euler-Lagrange extremals
to those which satisfy Weierstrass's necessary optimality condition \eqref{eqT:Weiers}).
In the present article we show that to formulate
Noether's theorem for admissible Lipschitz functions, one does not need to restrict
the set of Euler-Lagrange extremals to Pontryagin extremals, being
enough the restriction to those Euler-Lagrange extremals satisfying
the DuBois-Reymond condition:\footnote{In the autonomous case,
when the Lagrangian $L$ does not depend on the time variable
$t$ ($L = L(x,v)$), the DuBois-Reymond condition is known as the
\emph{second Erdmann equation}. For a survey of the classical
optimality conditions and related matters, we refer the reader
\textrm{e.g.} to \cite[Ch. 2]{MR91j:49001}.}
\begin{equation}
\label{eqT:DBR}
\frac{\partial L}{\partial t}\left(t,x(t),\dot{x}(t)\right) =
\frac{\mathrm{d}}{\mathrm{d}t} \left\{
L\left(t,x(t),\dot{x}(t)\right)
- \frac{\partial L}{\partial v}\left(t,x(t),\dot{x}(t)\right) \cdot \dot{x}(t)
\right\} \, .
\end{equation}
We remark that the DuBois-Reymond first-order necessary optimality
condition \eqref{eqT:DBR} is valid when $\mathcal{X}$ is the class
of Lipschitz functions,
and that \eqref{eqT:DBR} is a consequence of the Euler-Lagrange
and Weierstrass conditions.
For $\mathcal{X} = C^2$, \eqref{eqT:DBR}
follows from the Euler-Lagrange equation \eqref{eqT:EL} alone
(every Euler-Lagrange $C^2$-extremal is a DuBois-Reymond
$C^2$-extremal),\footnote{The opposite is not in general true.
For an example of a DuBois-Reymond $C^2$-extremal
which is not an Euler-Lagrange extremal see \cite[p.~24]{MR98b:49002a}.}
and therefore our generalization of Noether's theorem to the class
of Lipschitz functions (\textrm{cf.} \S\ref{SecT:MR}) gives
a proper extension of the smooth result.


\section{Review of Noether's Symmetry Theorem}

The universal principle described by Noether,
asserts that the invariance of a
problem with respect to a one-parameter group of transformations
implies the existence of a conserved quantity along
the smooth Euler-Lagrange extremals.

\begin{definition}[\textrm{cf.} \cite{MR2000m:49002}]
\label{Torres:defClassInv1}
If $C^2 \ni h^s(t,x) = \left(h_t^s(t),h_x^s(x)\right) : [a,b] \times \mathbb{R}^n
\rightarrow \mathbb{R} \times \mathbb{R}^n$,
$s \in \left(-\varepsilon,\varepsilon\right)$;
$h^0(t,x) = (t,x)$ for all $(t,x) \in [a,b] \times \mathbb{R}^n$; and
\begin{equation}
\label{eq:InvLag}
\int_{h_t^0(a)}^{h_t^0(b)}
L\left(t^s,h_x^s\left(x(t^s)\right),
\frac{\mathrm{d}}{\mathrm{d}t^s}h_x^s\left(x(t^s)\right)\right)\mathrm{d}t^s
= \int_a^b L\left(t,x(t),\dot{x}(t)\right)\mathrm{d}t \, ,
\end{equation}
for $t^s = h_t^s(t)$, all $s \in (-\varepsilon,\varepsilon)$, and all
$x(\cdot) \in C^2\left([a,b];\mathbb{R}^n\right)$;
then the fundamental problem of the calculus of variations is said to be invariant
under $h^s$.
\end{definition}

\begin{theorem}[Classical Noether's Theorem]
\label{Torres:ClassFNT}
If the fundamental problem of the calculus of variations is invariant under
$h^s$, in the sense of  Definition~\ref{Torres:defClassInv1}, then
\begin{multline}
\label{Torres:eq:equivClassFNT}
\frac{\partial L}{\partial v}\left(t,x(t),\dot{x}(t)\right) \cdot
\frac{\partial}{\partial s}\left.h_x^s\left(x(t)\right)\right|_{s = 0} \\
+ \left[L\left(t,x(t),\dot{x}(t)\right)
- \frac{\partial L}{\partial v}\left(t,x(t),\dot{x}(t)\right) \cdot \dot{x}(t)\right]
\frac{\partial}{\partial s} \left.h_t^s\left(t\right)\right|_{s = 0}
\end{multline}
is constant in $t \in [a,b]$ along every Euler-Lagrange $C^2$-extremal.
\end{theorem}

The following result for nonsmooth extremals,
is a trivial corollary from the optimal control results
in \cite{MR1901565,delfimEJC}.

\begin{theorem}
\label{Torres:COtoCV}
If the fundamental problem of the calculus of variations is invariant under
$h^s$, in the sense that Definition~\ref{Torres:defClassInv1} holds with
\eqref{eq:InvLag} satisfied for all $x(\cdot) \in Lip\left([a,b];\mathbb{R}^n\right)$,
then \eqref{Torres:eq:equivClassFNT}
is constant in $t$ along every $x(\cdot) \in Lip\left([a,b];\mathbb{R}^n\right)$
satisfying simultaneously the Euler-Lagrange \eqref{eqT:EL}
and Weierstrass \eqref{eqT:Weiers} necessary conditions.
\end{theorem}

Theorem~\ref{Torres:COtoCV} restrict the conclusion of
Noether's theorem to Pontryagin extremals.
The following question comes immediately to mind:
\emph{Is it really necessary to restrict the set
of nonsmooth Euler-Lagrange extremals in order to guarantee
that \eqref{Torres:eq:equivClassFNT} is conserved?}
In Section~\ref{Sec:NELEMFTSNCQ}
we show that a restriction is indeed necessary:
we provide an example of a Lipschitz Euler-Lagrange extremal
which is not a Weierstrass extremal,
and which fails to preserve \eqref{Torres:eq:equivClassFNT}.

While Pontryagin extremals are a natural
choice in optimal control, in the context of the calculus
of variations such restriction seems to be unnatural:
Theorem~\ref{Torres:COtoCV} does not simplify
to Theorem~\ref{Torres:ClassFNT} in the $C^2$ smooth case
(Euler-Lagrange equation differs from Weierstrass's necessary condition
in the $C^2$ smooth case). This means that Theorem~\ref{Torres:COtoCV}
does not give a proper extension of the classical Noether's theorem.
In Section~\ref{SecT:MR} we give a proper restriction
of the set of nonsmooth Euler-Lagrange extremals for which
Noether's theorem can still be asserted (Theorem~\ref{Th:MainResult}).


\section{Nonsmooth Euler-Lagrange Extremals may fail to satisfy
Noether's Conserved Quantities}
\label{Sec:NELEMFTSNCQ}

Let us consider the fundamental problem of the calculus of variations
with the Lagrangian given by $L(v) = (v^2 - 1)^2$: to minimize the
functional
\begin{equation*}
J\left[x(\cdot)\right] =
\int_0^1 \left\{\left[\dot{x}(t)\right]^2 - 1\right\}^2 \mathrm{d}t
\end{equation*}
over the class $\mathcal{X} = Lip$ of Lipschitz functions $x(\cdot)$
on the interval $[0,1]$ satisfying $x(0) = x(1) = 0$.
From the Euler-Lagrange equation \eqref{eqT:EL} one obtains that
any solution of this problem must satisfy
\begin{equation}
\label{eqT:EL:ex1}
\left\{\left[\dot{x}(t)\right]^2 - 1\right\} \dot{x}(t) = \text{const}
\quad a.e. \text{ on } [0,1] \, .
\end{equation}
As far as the problem is time-invariant (one can choose
$h_t^s(t) = t+s$, $h_x^s(x) = x$ in Definition~\ref{Torres:defClassInv1}),
Noether's conserved quantity \eqref{Torres:eq:equivClassFNT}
coincides with the DuBois-Reymond condition \eqref{eqT:DBR}:
\begin{equation*}
\left\{\left[\dot{x}(t)\right]^2 - 1\right\}^2
- 4 \left\{\left[\dot{x}(t)\right]^2 - 1\right\} \left[\dot{x}(t)\right]^2
= \text{const} \quad a.e. \text{ on } [0,1] \, ,
\end{equation*}
that is,
\begin{equation}
\label{eqT:DBR:ex1}
\left\{\left[\dot{x}(t)\right]^2 - 1\right\}
\left\{1 + 3\left[\dot{x}(t)\right]^2\right\} = \text{const}
\quad a.e. \text{ on } [0,1] \, .
\end{equation}
As is easily seen, there exist Lipschitz solutions of \eqref{eqT:EL:ex1}
which are not solutions of \eqref{eqT:DBR:ex1} (there are Lipschitz Euler-Lagrange
extremals which are not DuBois-Reymond extremals, and which not preserve the
quantity \eqref{Torres:eq:equivClassFNT} of Noether's theorem).
In fact, any Lipschitz function $x(\cdot)$ satisfying
$\dot{x}(t) \in \left\{-1,0,1\right\}$, $t \in [0,1]$, is an Euler-Lagrange
extremal. Among them, only $\dot{x}(t) \equiv 0$ and those with $\dot{x}(t)\pm 1$
satisfy \eqref{eqT:DBR:ex1}.


\section{Main Result}
\label{SecT:MR}

We formulate our Noether theorem for nonsmooth extremals
under a more general notion of invariance than the one in
Definition~\ref{Torres:defClassInv1}. We require the symmetry
transformation to leave the problem invariant up to first
order terms in the parameter, and up to exact
differentials.\footnote{The exact differentials are called
\emph{gauge-terms} in the literature (\textrm{cf.} \cite{MR83c:70020}).}

\begin{definition}
\label{def:QIUGT}
The integral functional \eqref{eqT:J} is quasi-invariant
under a one-parameter group of $C^1$-transformations
$(t,x) \longrightarrow \left(T(t,x,\dot{x},s),X(t,x,\dot{x},s)\right)$,
$|s| < \varepsilon$,
up to the gauge-term $\Phi(t,x,\dot{x})$ if, and only if,
\begin{multline}
\label{eqT:defQI}
\frac{\mathrm{d}}{\mathrm{d}t} \Phi\left(t,x(t),\dot{x}(t)\right)
= \frac{\mathrm{d}}{\mathrm{d}s} \Biggl\{
L\Biggl(T(t,x(t),\dot{x}(t),s),X(t,x(t),\dot{x}(t),s), \\
\left.\left.\left.
\frac{\frac{\mathrm{d}X(t,x(t),\dot{x}(t),s)}{\mathrm{d}t}}{\frac{\mathrm{d}
T(t,x(t),\dot{x}(t),s)}{\mathrm{d}t}}\right)\frac{\mathrm{d}
T\left(t,x(t),\dot{x}(t),s\right)}{\mathrm{d}t}\right\}\right|_{s = 0}
\end{multline}
for all $x(\cdot) \in Lip\left([a,b]; \mathbb{R}^n\right)$.
\end{definition}

\begin{remark}
As in the classical context,
we are assuming that the parameter transformations
$(t,x) \longrightarrow \left(T(t,x,\dot{x},s),X(t,x,\dot{x},s)\right)$
reduce to the identity for $s = 0$, that is,
\begin{equation*}
T(t,x,v,0) = t \, , \quad X(t,x,v,0) = x \, ,
\end{equation*}
for any choice of $t$, $x$, and $v$.
\end{remark}

\begin{remark}
It is obvious that the invariance notion
used in connection with Theorem~\ref{Torres:COtoCV}
implies the quasi-invariance up to a gauge-term
in Definition~\ref{def:QIUGT}.
\end{remark}

\begin{remark}
In the 1918 original paper of Emmy Noether \cite{JFM46.0770.01,MR53:10538},
Noether explains that the derivatives
of the trajectories $x$ may also occur in the
parameter group of transformations.
This possibility has been widely forgotten in the literature
of the calculus of variations,
the only exception seeming to be the textbook of
I. M. Gelfand and S. V. Fomin \cite{MR28:3353}.
Such possibility is, however, very interesting
from the point of view of optimal control
(\textrm{cf.} \cite{MR1901565,delfimEJC})
and is included in Definition~\ref{def:QIUGT}.
An example with relevance in Physics,
showing that the dependence of the invariance
transformations on the derivatives can be crucial in order to
obtain a conservation law, can be found in \cite{MoyoLeach2002}.
\end{remark}

\begin{theorem}[Noether's theorem for Lipschitz functions]
\label{Th:MainResult}
If \eqref{eqT:J} is quasi-invariant under the one-parameter
group of time-space transformations
\begin{equation*}
(t,x) \longrightarrow \left(T(t,x,\dot{x},s),X(t,x,\dot{x},s)\right)
\end{equation*}
up to the gauge-term $\Phi\left(t,x,\dot{x}\right)$, then
\begin{multline*}
\left[
L\left(t,x(t),\dot{x}(t)\right)
- \frac{\partial L}{\partial v}\left(t,x(t),\dot{x}(t)\right) \cdot \dot{x}(t)
\right]
\left.\frac{\partial}{\partial s} T(t,x(t),\dot{x}(t),s)\right|_{s = 0} \\
+ \frac{\partial L}{\partial v}\left(t,x(t),\dot{x}(t)\right)
\cdot \left.\frac{\partial}{\partial s} X(t,x(t),\dot{x}(t),s)\right|_{s = 0}
- \Phi(t,x(t),\dot{x}(t))
\end{multline*}
is constant in $t \in [a,b]$ along any
$x(\cdot) \in Lip\left([a,b]; \mathbb{R}^n\right)$ satisfying
\eqref{eqT:EL} and \eqref{eqT:DBR} (along any Lipschitz
Euler-Lagrange extremal which is also a Lipschitz DuBois-Reymond extremal).
\end{theorem}

\begin{proof}
Having in mind that for $s = 0$ we have the identity transformation,
$T(t,x,\dot{x},0) = t$, $X(t,x,\dot{x},0) = x$, condition
\eqref{eqT:defQI} yields
\begin{multline}
\label{eqT:diffDefQI}
\frac{\mathrm{d}}{\mathrm{d}t} \Phi\left(t,x(t),\dot{x}(t)\right) =
\frac{\partial L}{\partial t}\left(t,x(t),\dot{x}(t)\right)
\left.\frac{\partial}{\partial s} T\left(t,x(t),\dot{x}(t),s\right)\right|_{s = 0} \\
+ \frac{\partial L}{\partial x}\left(t,x(t),\dot{x}(t)\right) \cdot
\left.\frac{\partial}{\partial s} X\left(t,x(t),\dot{x}(t),s\right)\right|_{s = 0} \\
+ \frac{\partial L}{\partial v}\left(t,x(t),\dot{x}(t)\right) \cdot
\left( \frac{\mathrm{d}}{\mathrm{d}t}
\left.\frac{\partial}{\partial s}
X\left(t,x(t),\dot{x}(t),s\right)\right|_{s = 0} \right.\\
\left. - \dot{x}(t) \cdot \frac{\mathrm{d}}{\mathrm{d}t}
\left.\frac{\partial}{\partial s}
T\left(t,x(t),\dot{x}(t),s\right)\right|_{s = 0}\right) \\
+ L\left(t,x(t),\dot{x}(t)\right) \frac{\mathrm{d}}{\mathrm{d}t}
\left.\frac{\partial}{\partial s} T\left(t,x(t),\dot{x}(t),s\right)\right|_{s = 0}\, .
\end{multline}
From \eqref{eqT:EL} one can write
\begin{multline}
\label{eqT:FromEL}
\frac{\partial L}{\partial x}\left(t,x(t),\dot{x}(t)\right) \cdot
\left.\frac{\partial}{\partial s} X\left(t,x(t),\dot{x}(t),s\right)\right|_{s = 0} \\
+ \frac{\partial L}{\partial v}\left(t,x(t),\dot{x}(t)\right) \cdot
\frac{\mathrm{d}}{\mathrm{d}t}
\left.\frac{\partial}{\partial s} X\left(t,x(t),\dot{x}(t),s\right)\right|_{s = 0} \\
= \frac{\mathrm{d}}{\mathrm{d}t} \left(
\frac{\partial L}{\partial v}\left(t,x(t),\dot{x}(t)\right) \cdot
\left.\frac{\partial}{\partial s} X\left(t,x(t),\dot{x}(t),s\right)\right|_{s = 0}
\right) \, ,
\end{multline}
while from \eqref{eqT:DBR} one gets
\begin{multline}
\label{eqT:FromDBR}
\frac{\partial L}{\partial t}\left(t,x(t),\dot{x}(t)\right) \cdot
\left.\frac{\partial}{\partial s} T\left(t,x(t),\dot{x}(t),s\right)\right|_{s = 0}\\
+ L\left(t,x(t),\dot{x}(t)\right) \frac{\mathrm{d}}{\mathrm{d}t}
\left.\frac{\partial}{\partial s} T\left(t,x(t),\dot{x}(t),s\right)\right|_{s = 0}\\
- \frac{\partial L}{\partial v}\left(t,x(t),\dot{x}(t)\right) \cdot \dot{x}(t)
\frac{\mathrm{d}}{\mathrm{d}t}
\left.\frac{\partial}{\partial s} T\left(t,x(t),\dot{x}(t),s\right)\right|_{s = 0}\\
= \frac{\mathrm{d}}{\mathrm{d}t} \left\{\left(L\left(t,x(t),\dot{x}(t)\right)
- \frac{\partial L}{\partial v}\left(t,x(t),\dot{x}(t)\right) \cdot \dot{x}(t)\right)
\left.\frac{\partial}{\partial s} T\left(t,x(t),\dot{x}(t),s\right)\right|_{s = 0}
\right\} \, .
\end{multline}
Substituting \eqref{eqT:FromEL} and \eqref{eqT:FromDBR} into \eqref{eqT:diffDefQI},
\begin{multline*}
\frac{\mathrm{d}}{\mathrm{d}t} \left\{
\frac{\partial L}{\partial v}\left(t,x(t),\dot{x}(t)\right)
\cdot \left.\frac{\partial}{\partial s} X(t,x(t),\dot{x}(t),s)\right|_{s = 0}
- \Phi(t,x(t),\dot{x}(t)) \right. \\
\left. + \left(
L\left(t,x(t),\dot{x}(t)\right)
- \frac{\partial L}{\partial v}\left(t,x(t),\dot{x}(t)\right) \cdot \dot{x}(t)\right)
\left.\frac{\partial}{\partial s} T(t,x(t),\dot{x}(t),s)\right|_{s = 0} \right\} = 0 \, ,
\end{multline*}
and the pretended conclusion is obtained.
\end{proof}



\end{document}